\documentclass[12pt,twoside]{article}
\usepackage[english]{babel}
\usepackage[latin1]{inputenc}
\usepackage{amsmath}
\usepackage{cite}
\usepackage{amssymb,amsfonts}
\usepackage{graphicx}
\usepackage{times,amssymb,amscd}

\newtheorem{theorem}{Theorem}[section]

\newtheorem{exmple}[theorem]{Example}
\newtheorem{defn}[theorem]{Definition}
\newtheorem{rmrk}[theorem]{Remark}

\newcommand{\bA}{\mathbf{A}}

\newcommand{\bH}{\mathbf{H}}

\newcommand{\bL}{\mathbf{L}}
\newcommand{\bN}{\mathbf{N}}
\newcommand{\bR}{\mathbf{R}}
\newcommand{\bS}{\mathbf{S}}
\newcommand{\bV}{\mathbf{V}}

\newcommand{\bX}{\mathbf{X}}

\newcommand{\BV}{\boldsymbol{V}}

\newcommand{\Ba}{\boldsymbol{a}}
\newcommand{\Bb}{\boldsymbol{b}}

\newcommand{\Bh}{\boldsymbol{h}}
\newcommand{\Bu}{\boldsymbol{u}}
\newcommand{\Bv}{\boldsymbol{v}}

\newcommand{\BF}{\boldsymbol{F}}
\newcommand{\BQ}{\boldsymbol{Q}}

\newcommand{\BB}{\boldsymbol{B}}
\newcommand{\BE}{\boldsymbol{E}}
\newcommand{\BCw}{\boldsymbol{Cw}}
\newcommand{\BW}{\boldsymbol{W}}
\newcommand{\BA}{\boldsymbol{A}}

\newcommand{\BH}{\boldsymbol{H}}
\newcommand{\BX}{\boldsymbol{X}}
\newcommand{\BY}{\boldsymbol{Y}}
\newcommand{\BU}{\boldsymbol{U}}
\newcommand{\Bs}{\boldsymbol{s}}
\newcommand{\Bc}{\boldsymbol{c}}

\newcommand{\cP}{\mathcal{P}}

\newcommand{\cM}{\mathcal{M}}

\newcommand{\HYP}{\bH^3}

\newcommand{\SLR}{\widetilde{\bS\bL_2\bR}}
\newcommand{\NIL}{\mathbf{Nil}}
\newcommand{\SOL}{\mathbf{Sol}}

\begin{document}
\pagestyle{myheadings}
\markboth{\centerline{Emil Moln\'ar and Jen\H o Szirmai}}
{Dense ball packings by tube manifolds}
\title
{Dense ball packings by tube manifolds as new
models for hyperbolic
crystallography \\ \vspace{5mm} \normalsize{\it To the 200th Anniversary of János Bolyai's
Absolute Geometry}
}
\vspace{5mm}
\author{Emil Moln\'ar and Jen\H o Szirmai \\
\normalsize Department of Algebra and Geometry, Institute of Mathematics,\\
\normalsize Budapest University of Technology and Economics, \\
\normalsize M\H uegyetem rkp. 3., H-1111 Budapest, Hungary \\
\normalsize emolnar@math.bme.hu,~szirmai@math.bme.hu
\date{\normalsize{\today}}}

\maketitle
\begin{abstract}
We intend to continue our previous papers (\cite{MSz17} and \cite{MSz18}, as indicated there)
on dense ball packing hyperbolic space $\HYP$ by equal balls, but here with
centres belonging to different orbits of the fundamental group
$Cw(2z, 3 \le z \in \bN$, odd number),
of our new series of {\it tube or cobweb
manifolds} $Cw = \HYP/\BCw$ with $z$-rotational symmetry. As we know,
$\BCw$ is a fixed-point-free isometry group, acting on $\HYP$ discontinuously
with appropriate tricky fundamental domain $Cw$, so that every point has
a ball-like neighbourhood in the usual factor-topology.

Our every $Cw(2z)$ is minimal, i.e. does not cover regularly a
smaller manifold. It can be derived by its general symmetry group
$\BW(u; v; w = u)$ that is a complete Coxeter orthoscheme reflection group,
extended by the half-turn $\Bh$ $(0 \leftrightarrow 3, 1 \leftrightarrow 2)$
of the complete orthoscheme
$A_0A_1A_2A_3 \sim  b_0b_1b_2b_3$ (Fig.~1). The vertices $A_0$ and $A_3$ are outer
points of
the (Beltrami-Cayley-Klein) {\it B-C-K model of} $\HYP$, as
$1/u + 1/v \le 1/2$ is required, $3 \le  u = w, v$ for the above orthoscheme
parameters. For the above simple manifold-construction
we specify $u = v = w = 2z$. Then the polar planes
$a_0$ and $a_3$ of $A_0$ and $A_3$, respectively, make complete with reflections
$\Ba_0$ and $\Ba_3$ the Coxeter reflection group, where the other
reflections
are denoted by $\Bb^0$, $\Bb^1$, $\Bb^2$, $\Bb^3$ in the sides of the
orthoscheme $b^0b^1b^2b^3$.

The situation is described first in Figure 1 of the half trunc-orthoscheme
$W$ and its usual  extended Coxeter diagram, moreover, by
the scalar product matrix $(b^{ij}) = (\langle \Bb^i, \Bb^j \rangle)$ in
formula (1.1) and its inverse
$(A_{jk}) = (\langle \BA_j, \BA_k \rangle)$ in (1.3).
These will describe the hyperbolic angle and
distance metric of the half trunc-orthoscheme $W$, then its ball packings,
densities, then those of the manifolds $Cw(2z)$.

As first results we concentrate only on particular constructions by computer for probable material model
realizations, atoms or molekules by equal balls, for general $W(u;v;w=u)$ as well, summarized at the end of our paper.

\end{abstract}
%\tableofcontents
%
\section{Introduction}
\subsection{Hyperbolic space $\HYP$, as a projective metric space. Generalized Beltrami-Cayley-Klein
(B-C-K) model}

Our Bolyai-Lobachevsky hyperbolic space $\HYP =\cP^3\cM(\bV^4,\BV_4,\bR,\sim^,
\langle~,~\rangle)$
will be a projective-metric space over a real vector space $\bV^4$ for points
$X (\BX \sim c\BX,~c \in \bR \setminus\{0\})$; its dual (i.e. linear form space)
$\BV_4$
will describe planes ($2$-planes) $u (\Bu \sim \Bu c,~ c \in \bR \setminus \{0\})$.
The scalar product $\langle~,~\rangle$ will be specified in (1.1), (1.3) by a
so-called complete orthoscheme as a projective coordinate
simplex  $b^0b^1b^2b^3 = A_0 A_1 A_2 A_3$   by $\bA_ib^j = \delta_i^j$
(Kronecker delta), i.e. $b^i = A_jA_kA_l, ~ \{i; j; k; l\} = \{0; 1; 2; 3\}$.

Furthermore, the starting Coxeter-Schl\"afli matrix (1.1) will be defined by three natural parameters
$3 \le u, v, w;$ (think of $u=5, v=3, w=5$ at the characteristic orthoscheme of our previous football manifold
$\{5,6,6\}$ in \cite{MSz17}, \cite{MSz18});

\begin{equation}
\begin{gathered}
(b^{ij})
:=(\langle \Bb^i,\Bb^j \rangle )=(\cos(\pi-\beta^{ij})),~ \beta^{ii}=\pi;\\
(b^{ij})=\begin{pmatrix}
1 &-\cos\frac{\pi}{u} & 0&0 \\
-\cos \frac{\pi}{u} &1& -\cos\frac{\pi}{v}&0 \\
0 & -\cos\frac{\pi}{v} &1&-\cos \frac{\pi}{w} \\
0&0&-\cos \frac{\pi}{w} &1
\end{pmatrix},
 \end{gathered} \tag{1.1}
\end{equation}
as scalar products of basis forms $\Bb^i \in \BV_4$ $(i = 0, 1, 2, 3)$ to
the side faces of the coordinate simplex $b^0b^1b^2b^3 = A_0 A_1 A_2 A_3$ as usual.
Thus, the essential face angles $(\angle b^ib^j)=\beta^{ij}$ are $\beta^{01} = \pi/u = \pi/5$,
$\beta^{12} = \pi/v = \pi/3$, $\beta^{23} = \pi/w = \pi/5$, the others are
$\beta^{02} = \beta^{03} = \beta^{13} = \pi/2$ (rectangle).

The above scalar product by (1.1) is equivalent by giving a (polar) plane $\rightarrow$ (pole) point
$\BV_4 \rightarrow \bV^4$ (linear symmetric) polarity
\begin{equation}
\begin{gathered}
b \rightarrow B,~\Bb^i \rightarrow \BB^i:=b^{ij}\bA_j,~i,j\in\{0;1;2;3\},~
\text{Einstein-Schouten convention)};\\
\text{i.e. by (1.1)}: ~\Bb^0 \rightarrow \BB^0=\BA_0-\cos\frac{\pi}{u}\BA_1,~\tag{1.2}\\
\Bb^1 \rightarrow \BB^1=-\cos\frac{\pi}{u}\BA_0+\BA_1-\cos\frac{\pi}{v}\BA_2, \\
\Bb^2 \rightarrow \BB^2=-\cos\frac{\pi}{v}\BA_1+\BA_2-\cos\frac{\pi}{w}\BA_3,~
\Bb^3 \rightarrow \BB^3=-\cos\frac{\pi}{w}\BA_2+\BA_3.
\end{gathered}
\end{equation}
First for $\Bu = \Bb^i u_i$, $\Bv = \Bb^j v_j$ let $\langle \Bu, \Bv \rangle := (u_i\BB^i, \Bb^j v_j)$
$= (u_i b^{ik} \BA_k, \Bb^j v_j)$ $= u_i b^{ik} \delta_k^j v_j = u_i b^{ij} v_j$
be defined step by step, Then
$$
\cos(\angle uv) = \frac{-\langle \Bu, \Bv \rangle}{\sqrt{\langle \Bu, \Bu \rangle \langle \Bv, \Bv \rangle}}
$$
defines the ``usual angle" according to (1.1), (1.2), indeed.
Proper plane $u$ means that $\langle \Bu, \Bu \rangle = (\BU \Bu) > 0$.
Then its pole $U(\BU)$ is (unproper) outer point in (B-C-K) model. etc.

Its ``inverse scalar product" is by
\[
\begin{gathered}
(A_{ij})=(b^{ij})^{-1}=\langle \BA_i, \BA_j \rangle:=\\
=\frac{1}{B} \begin{pmatrix}
\sin^2{\frac{\pi}{w}}-\cos^2{\frac{\pi}{v}}& \cos{\frac{\pi}{u}}\sin^2{\frac{\pi}{w}}& \cos{\frac{\pi}{u}}\cos{\frac{\pi}{v}} & \cos{\frac{\pi}{u}}\cos{\frac{\pi}{v}}\cos{\frac{\pi}{w}} \\
\cos{\frac{\pi}{u}}\sin^2{\frac{\pi}{w}} & \sin^2{\frac{\pi}{w}} & \cos{\frac{\pi}{v}}& \cos{\frac{\pi}{w}}\cos{\frac{\pi}{v}} \\
\cos{\frac{\pi}{u}}\cos{\frac{\pi}{v}} & \cos{\frac{\pi}{v}} & \sin^2{\frac{\pi}{u}}  & \cos{\frac{\pi}{w}}\sin^2{\frac{\pi}{u}}  \\
\cos{\frac{\pi}{u}}\cos{\frac{\pi}{v}}\cos{\frac{\pi}{w}}  & \cos{\frac{\pi}{w}}\cos{\frac{\pi}{v}} & \cos{\frac{\pi}{w}}\sin^2{\frac{\pi}{u}}  & \sin^2{\frac{\pi}{u}}-\cos^2{\frac{\pi}{v}}
\end{pmatrix}, \tag{1.3}
\end{gathered}
\]
where
$$
B=\det(b^{ij})=\sin^2{\frac{\pi}{u}}\sin^2{\frac{\pi}{w}}-\cos^2{\frac{\pi}{v}} <0, \ \ \text{i.e.} \ \sin{\frac{\pi}{u}}\sin{\frac{\pi}{w}}-\cos{\frac{\pi}{v}}<0.
$$
We assume, and this will be crucial in the following, that $u = w$, so our
orthoscheme will be symmetric by a half-turn
$\Bh$: $0 \leftrightarrow 3, 1\leftrightarrow 2$.
The half-turn axis $h$ joins
the midpoints $F_{03}$ of $A_0$, $A_3$ and $F_{12}$ of $A_1$, $A_2$
(see also Fig.~1).

Fig.~1 shows our novelty at later cobweb (tube) manifolds: $\pi/u +\pi/v < \pi/2$.
Thus, our scalar product (1.1) (and its inverse (1.3) as well) will be of
signature $(+ + + -)$, $A_0$ and $A_3$ will be outer vertices of the
B-C-K-model, so we truncate the orthoscheme by their proper polar planes $a_0$
and $a_3$ (look also at the extended Coxeter-Schläfli diagram in Fig.~1),
respectively, to obtain a compact domain with proper (or interior) points,
a so-called complete orthoscheme (trunc-orthoscheme). Algebraically,
the upper minor determinant sequence in (1.1) guarantees the signature
$(+ + + -)$
of the scalar products (1.1), (1.3) to be hyperbolic indeed.

Extesions of angle and distance metrics are also standard by complex $\cos$ and $\cosh$ functions
(through exponential one, of course), respectively: $\cosh x = \cos(x/i)$, $i$
is the imaginary unit. We mention only that the matrix (1.3) defines the scalar
product of basis vectors  of $\bV^4$. These determine the $XY$ distance in $\HYP$
through the scalar product $\langle \BX, \BY \rangle = X^i A_{ij} Y^j$
of vectors $\BX = X^i \BA_i$ and $\BY = Y^j \BA_j \in \bV^4$
(Einstein-Schouten index conventions). Namely,
\begin{equation}
\cosh{\frac{XY}{k}}=\frac{-\langle \BX,\BY \rangle}{\sqrt{\langle
\BX,\BX \rangle \langle \BY,\BY \rangle}}, ~(\langle \BX,\BX \rangle,~\langle
\BY,\BY \rangle <0). \tag{1.4}
\end{equation}
Here $k=\sqrt{-\frac{1}{K}}$ is the natural length unit, where $K$
is the sectional curvature. $K$ can be chosen to
$-1$, so $k = 1$ in the following. But other $K$ can be important in the applications (in nano size)!

We recall the volume formula of complete orthoscheme $\mathcal{O}(\beta^{01}, \beta^{12}, \beta^{23})$ by R.~Kellerhals to the
Coxeter-Schl\"afli matrix (1.1) on the genial ideas of N.~I. Lobachevsky.
\begin{theorem}{\rm{(Kellerhals \cite{K89},~Lobachevsky)}} The volume of a three-dimensional hyperbolic
complete orthoscheme $\mathcal{O}(\beta^{01}, \beta^{12}, \beta^{23}) \subset \mathbf{H}^3$
is expressed with the essential
angles $\beta^{01}=\frac{\pi}{u}$, $\beta^{12}=\frac{\pi}{v}$, $\beta^{23}=\frac{\pi}{w}$, $(0 \le \alpha_{ij}
\le \frac{\pi}{2})$ in the following form:

\begin{align}
&\mathrm{Vol}(\mathcal{O})=\frac{1}{4} \{ \mathcal{L}(\beta^{01}+\theta)-
\mathcal{L}(\beta^{01}-\theta)+\mathcal{L}(\frac{\pi}{2}+\beta^{12}-\theta)+ \notag \\
&+\mathcal{L}(\frac{\pi}{2}-\beta^{12}-\theta)+\mathcal{L}(\beta^{23}+\theta)-
\mathcal{L}(\beta^{23}-\theta)+2\mathcal{L}(\frac{\pi}{2}-\theta) \}, \tag{1.5}
\end{align}
where $\theta \in [0,\frac{\pi}{2})$ is defined by:
$$
\tan(\theta)=\frac{\sqrt{ \cos^2{\beta^{12}}-\sin^2{\beta^{01}} \sin^2{\beta^{23}
}}} {\cos{\beta^{01}}\cos{\beta^{23}}},
$$
and where $\mathcal{L}(x):=-\int\limits_0^x \log \vert {2\sin{t}} \vert dt$ \ denotes the
Lobachevsky function (introduced by J.~Milnor in this form).
\end{theorem}
The volume $\mathrm{Vol}(B(R))$ of a {\it ball} $B(R)$ of radius $R$ can be computed by the classical formula of J.~Bolyai:
\begin{equation}
\begin{gathered}
\mathrm{Vol}(B(R))=2\pi (\cosh(R) \sinh(R) -R)=\pi(\sinh(2R)-2R)=\\
=\frac{4}{3} \pi R^3(1+\frac{1}{5}R^2+\frac{2}{105}R^4+\dots). \tag{1.6}
\end{gathered}
\end{equation}
The usual (plane) reflection $\sigma(u, U):
X(\bX) \rightarrow Y(\BY)$ in the plane $u(\Bu)$ with pole $U(\BU)$ is for points:
\begin{equation}
\sigma(u,U):\BX \rightarrow \BY=\BX - \frac{2(\BX\Bu)}{\langle \Bu,\Bu \rangle}
\BU,~\text{here} ~ \langle \Bu,\Bu \rangle=(\BU\Bu).  \tag{1.7}
\end{equation}
\begin{figure}[ht]
\begin{center}
\includegraphics[width=6cm]{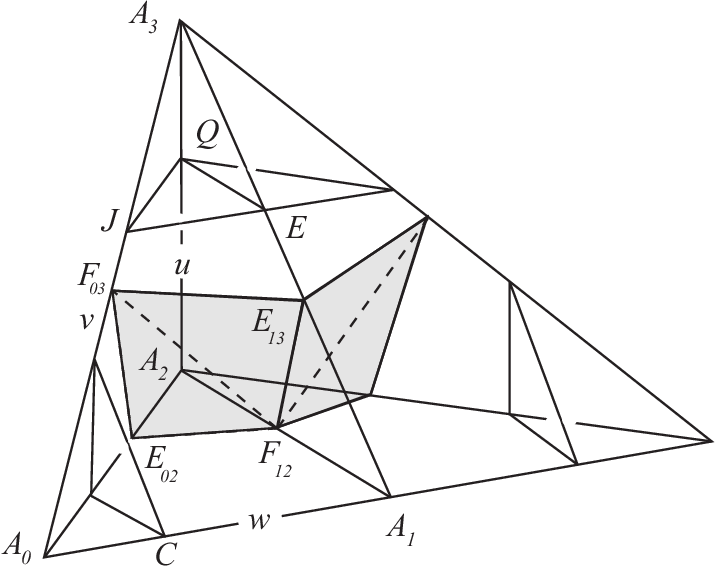} \includegraphics[width=7cm]{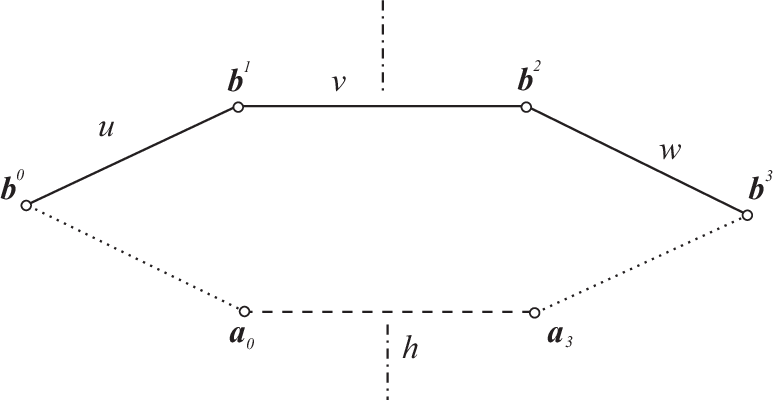}
\end{center}
\caption{The half trunc-orthoscheme fundamental domain $W(u; v; w=u)$, i.e. that of
the extended
(by halfturn $\Bh$)
complete orthoscheme group and its extended Coxeter-Schl\"afli diagram.}
\end{figure}
Thus we could easily describe the extended complete
orthoscheme reflection group $\BW(u; v; w=u)$.
A point $X(\BX)$ has an orthogonal projection $X_b(\BX_b)$ on a plane $b$
with pole $B$ by
\begin{equation}
\BX_b=\BX - \frac{(\BX\Bb)}{\langle \Bb,\Bb \rangle}\BB,~
\text{by} ~ \langle \Bb,\Bb \rangle=(\BB\Bb). \tag{1.8}
\end{equation}
So the distance of point $X$ from plane $b$, denoted simply by $Xb$,
can easily be expressed, etc.
\subsection{The cobweb (tube) manifolds $Cw(2z, 2z, 2z)$ $= Cw(2z)$}

For our new cobweb manifolds we start with its previously mentioned symmetry group
$\BW(u; v; w = u)$ as complete extended reflection group. We shall apply
$u = v = w = 2z$, $3 \le z$ is odd natural number,
i.e. $b^0b^1b^2b^3 = A_0 A_1 A_2 A_3$
is an orthoscheme; it will be complete,
i.e. doubly truncated with polar planes
$a_3$ and $a_0$ of $A_3$ and $A_0$, respectively, called also
half trunc-orthoscheme $W(2z)$.
\begin{figure}[]
\begin{center}
\includegraphics[width=10cm]{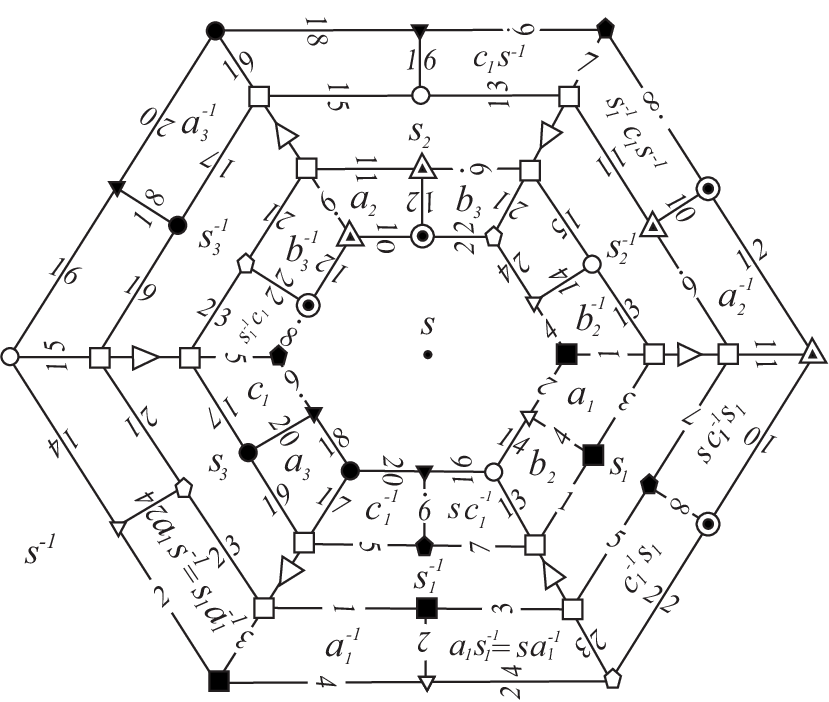}
\end{center}
\caption{The starting cobweb (tube) manifold $Cw(6)$ with its symbolic face pairing
isometries: e.g. $s^{-1} \rightarrow s$ by mapping (screw motion) $\Bs$. Signed equivalent edge triples are equally
numbered (from 1 to 24, odd and even edges play different roles).
We obtain signed vertex classes (by various symbols), indicated here
all together $1 + 3 \times 3 = 10$ ones. Any point has a
ball-like neighbourhood (for later nanotube manifold by the above fundamental group
$\BCw(6)$, on the base of $\BW(6), z = 3$.
}
\end{figure}
We consider this smaller asymmetric unit $W(2z=u=v=w)$
with the half-turn axis $h$
and a variable halving plane through $h$. Then we choose the point
$Q = a_3 \cap A_3A_0$ with
its stabilizer subgroup $\BW_Q$ of order $4u = 8z$ in the extended reflection group
$\BW(2z)$ to its fundamental domain $W(2z)$ (Fig.~1).
Then we reflect $W(2z)$ around $Q$ to get the cobweb polyhedron $Cw(2z)$
(Fig.~2 first for $z=3$)
as a fundamental domain of our new manifold with a new interesting fundamental
group denoted by $\BCw(2z)$.

Carbon atoms, or other ones (with valence $4$, e.g. at $F_{12}$
and its $\BCw$-equivalent positions) can be placed very naturally in this
tube-like structure.

It turns out that for $z = 4p-1$ and for $z = 4q+1$ $(1 \le p, q \in \bN )$
we get two analogous series (Fig.~6,7), each of them is unique by
$z$-rotational symmetry and manifold requirements.

Then come our new ball packing constructions as our new initiatives with more
ball centre orbits by the above symmetry groups $\BW(u; v; w = u)$
but equal balls.
There is only a part of a ball $B_X$ in half orthoschem $W$ of centre $X \in W$
just by the order $|\BW_X|$ of stabilizer subgroup $\BW_X$.
Our goal, to find the densest ball packing on a
natural average, is not completed yet. Our top density here is $0.68248 \dots$??
This problem seems to be very hard, of course (maybe hopeless, in general?).

Our summaries follow in Theorems 4.1, 4.2 in Section 4.
\begin{figure}[]
\begin{center}
\includegraphics[width=7cm]{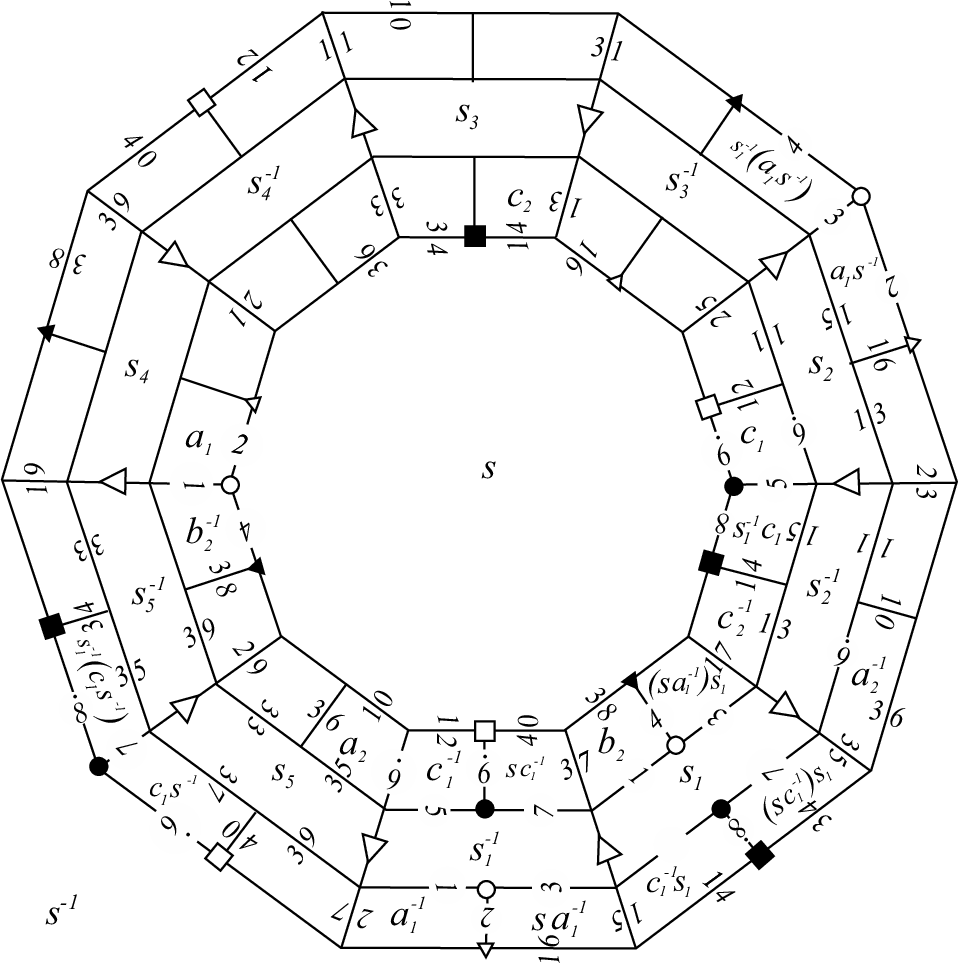} \includegraphics[width=5.5cm]{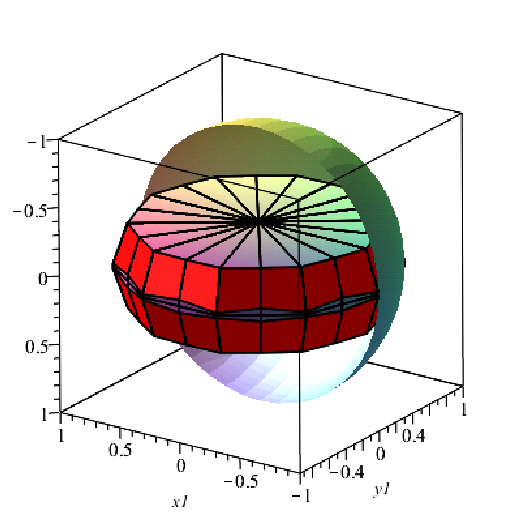}
\end{center}
\caption{Our old cobweb (tube) manifolds $Cw(2z = 8q+2)$, illustrated by $z = 5$, $q = 1$.
A picture of its animation in Beltrami-Cayley-Klein model.
}
\end{figure}
\section{Manifold constructions in new versions, with new ball packing initiatives}
\begin{figure}[]
\begin{center}
\includegraphics[width=8cm]{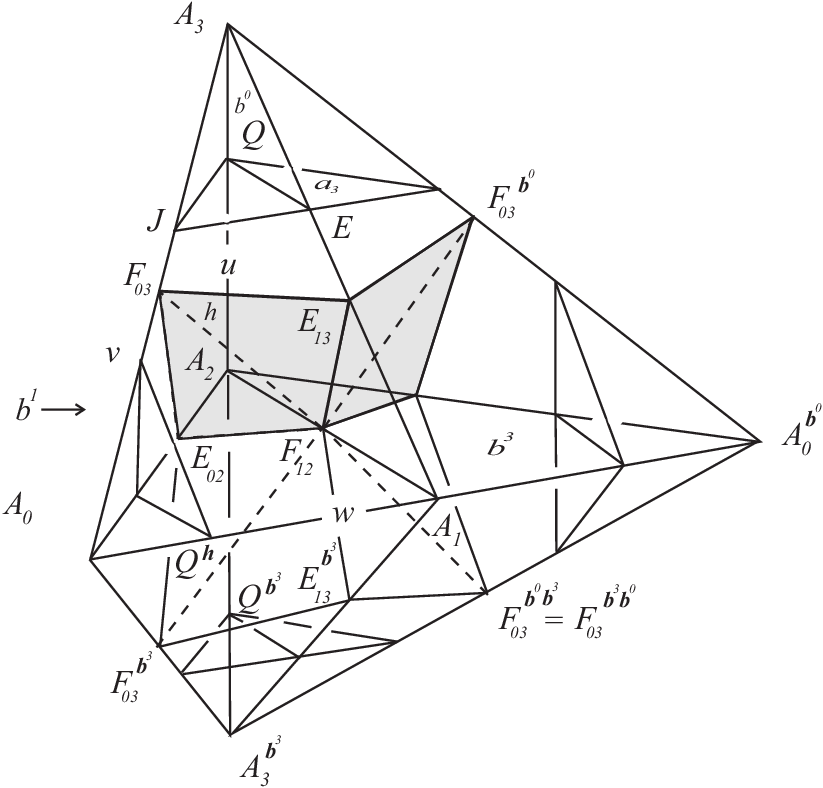}
\end{center}
\caption{The $\Bb^3$-reflected part $W(u;v;w=u)$ and its $\Bb^0$ mirror
image $W^{\Bb^0}$ in the coordinate plane $b^0$.
So the new bent ``quadrangle" face
$F_{03}E_{13}F_{03}^{\Bb^0}E_{02}^{\Bb^0}F_{03}^{\Bb^0 \Bb^3}E_{13}^{\Bb^3}F_{03}^{\Bb^3}E_{02}$
will be a $c_i$-type side face of the new ``cut-reglued" fundamental polyhedron
$Cw$. The $\Ba_3$ reflected part of $b^2$-face (together with $b^2$) will be an $s_j$-type face of the new $Cw$ (not indicated here).
}
\end{figure}
\subsection{Construction of cobweb (tube) manifold $Cw(6,6,6)=Cw(6)$ and the general two series $Cw(2z)$}
By the theory, e.g. in \cite{W} (cited also in works \cite{CaTe09}, \cite{P98},
\cite{S83}, we have to construct a fixed point free group acting in hyperbolic
space $\HYP$ with the above compact fundamental domain.
In the Introduction to Fig.~1 and analogously to Fig.~2 we have repeatedly described from
\cite{MSz18} the extended reflection group $\BW(6;6;6)=\BW(6)$ with fundamental
domain $W(6)$,
as a half complete Coxeter orthoscheme, and glued together to
the cobweb polyhedron $Cw(6;6;6)$=$Cw(6)$ as Dirichlet-Voronoi (in short $D-V$)
cell of the kernel point $Q$
by its orbit under the group $\BW(6)$. Now by Fig.~2,4,5 we shall give a new simpler
face identification of $Cw(6)$, so that it will be a cut-reglued
fundamental polyhedron
of the previous fixed-point-free group,
denoted also by $\BCw(6)$, generated just by the new face identifying
isometries (as hyperbolic screw motions).

By gluing $4u=24$ domains (at $A_2$ and) at $Q$ around (whose stabilizer subgroup $\BW_Q$ is just
of order $|\BW_Q|=4u=24$ as for $A_2$ as well), we simply "kill out" the fixed points of $\BW(6)$,
as we made in our former paper \cite{MSz18}.
\begin{figure}[]
\begin{center}
\includegraphics[width=11cm]{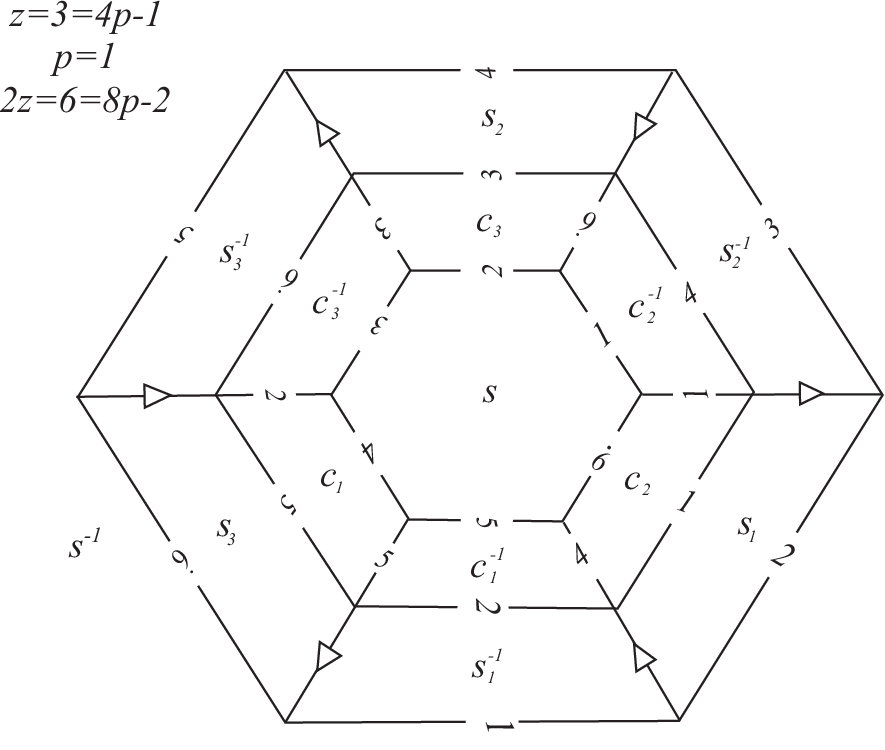}
\end{center}
\caption{The new simplified $Cw(6=2z)$ manifold. Only with $7=2z+1$ edge classes
($4$ edges in each of $6$ classes, $6=2z$ edges in class
$-\!-\!\!\vartriangleright$.) The former odd classes $1$, $3$ will
be the new class $1$ here; former
$5,7\rightarrow$ new edge class $2$. The former even classes are glued together in
the interior of $c_i$-type
faces. We have only $14=4z+2$ side faces, i.e. $7=2z+1$
generators for the fundamental group $\BCw$.
}
\end{figure}
We shortly repeat from \cite{MSz18}
our new unified algorithmic presentation of fundamental group of the first series.
At the beginning there stand the generator screw motions:
$\Bs$ for the tube form by rotational component $2\pi(z-1)/2z$; $\Bs_1,\dots, \Bs_z$
are half screws around; $\Bc_1, \dots, \Bc_z$ come from specific rotation around
the tube axis, combined with a halfturn (a conjugate of $\Bh$). Then come the relations
to the numbered or signed edges:
\begin{equation}
\begin{gathered}
\BCw(2z=8p-2):(\Bs,\Bs_1,\Bs_2,\dots,\Bs_z,\Bc_1,\Bc_2,\dots,\Bc_z;-\\
\text{edge}~2i-1:~\Bs_i \Bc^{-1}_{i+p}
\Bc^{-1}_{i-1+2p}\Bs^{-1}=1;~\text{edge}~2i:~\Bs^{-1}\Bc_i\Bc_{i-p}\Bs^{-1}=1;-\\
\text{edge}~-\!\!\!-\!\!\!\vartriangleright\!\!-\!\!-:
~1=\prod^i(\Bc_i\Bc_{i-p}\Bc_{i-1+2p}\Bc_{i+p}),~\text{for}~
i=1,2,\dots,z;~\text{indices are} \mod{z}).
\end{gathered} \tag{2.1}
\end{equation}
In case $p = 1$ we obtain the presentation of $\BCw(6)$ by Fig.~5.
We analogously obtain the algorithmic presentation of our second series
\begin{equation}
\begin{gathered}
\BCw(2z=8q+2):(\Bs,\Bs_1,\Bs_2,\dots,\Bs_z,\Bc,\Bc_1,\Bc_2,\dots,\Bc_z;-\\
\text{edge}~2i-1:~\Bs_i \Bc^{-1}_{i-q}
\Bc^{-1}_{i+2q}\Bs^{-1}=1;~\text{edge}~2i:~\Bs^{-1}\Bc_i\Bc_{i+q}\Bs^{-1}=1;-\\
\text{edge}~-\!\!\!-\!\!\!\vartriangleright\!\!\!-\!\!-:~1=
\prod^i(\Bc_i\Bc_{i+q}\Bc_{i+2q}\Bc_{i-q}),~\text{for}~ i=1,2,\dots,z;
~\text{indices are} \mod{z}).
\end{gathered} \tag{2.2}
\end{equation}
In case $q = 1$ we obtain the presentation of $\BCw(10)$ by Fig.~3, here in new interpretation.

In both cases we have some consequences, e.g. in first case
$\Bs_i = \Bs \Bc_{i-1+2p} \Bc_{i+p}  = \Bc_i \Bc_{i-p}  \Bs^{-1}$ that
was utilized in the last product relation.

The face pairing structure of these manifolds can be derived generally in Fig.~6,7
(see also Fig.~3).
The above tube screw motion $\Bs$
has rotation component $2\pi(z-1)/2z$, throughout in the following.
The crucial difference between the two series is that the third edge in class $1$
in the triple will be placed opposite to each other on the cobweb (tube) polyhedron,
as (2.1) and (2.2) express as well.
Geometrically each manifold $Cw(2z)$ and also $W(2z)$ appears in two mirror forms (in our figures as well), equivalent
mathematically, maybe not so in the occasional applications
(similarly as in the Euclidean crystallography).

These manifolds realize nanotubes in small (nanometer=$10^{-9}$m) size, very probably(!?).
And of course, there arise new open questions.

For completeness we recall the presentation of the symmetry group of our cobweb
(tube) manifolds $Cw(2z)$ in more general (a most economical) version.
Thus we have for the extended complete half-orthoscheme reflection group with parameters
$3 \le  u = w, v \in \bN$  by the Coxeter-Schl\"afli diagram in Fig.~1, as follows:

\begin{equation}
\begin{gathered}
\BW(u;v;w=u)=(\Ba_3,\Bb^0,\Bb^1, \Bh ~ -\!\!- 1=\Ba_3\Ba_3=\Bb^0\Bb^0=\Bb^1\Bb^1=\Bh\Bh=\\
=(\Ba_3\Bb^0)^2=(\Ba_3 \Bb^1)^2=(\Ba_3 \Bh \Bb^1 \Bh)^2=(\Bb^0\Bb^1)^u=\\
=(\Bh\Bb^0\Bh\Bb^0)^2=(\Bh\Bb^0\Bh\Bb^1)^2=
(\Bh\Bb^1\Bh\Bb^1)^v).
\end{gathered} \tag{2.3}
\end{equation}

Of course, we can express the generators of the above $\BCw(2z=8p-2)$ and $\BCw(2z=8q+2)$,
respectively, with generators of $\BW$, and check the corresponding relations
by those of $\BW$ in (2.3). Here the rotation $(\Bb^0\Bb^1)$ of
order $u (= 2z)$ plays important role.

\begin{figure}[]
\begin{center}
\includegraphics[width=9cm]{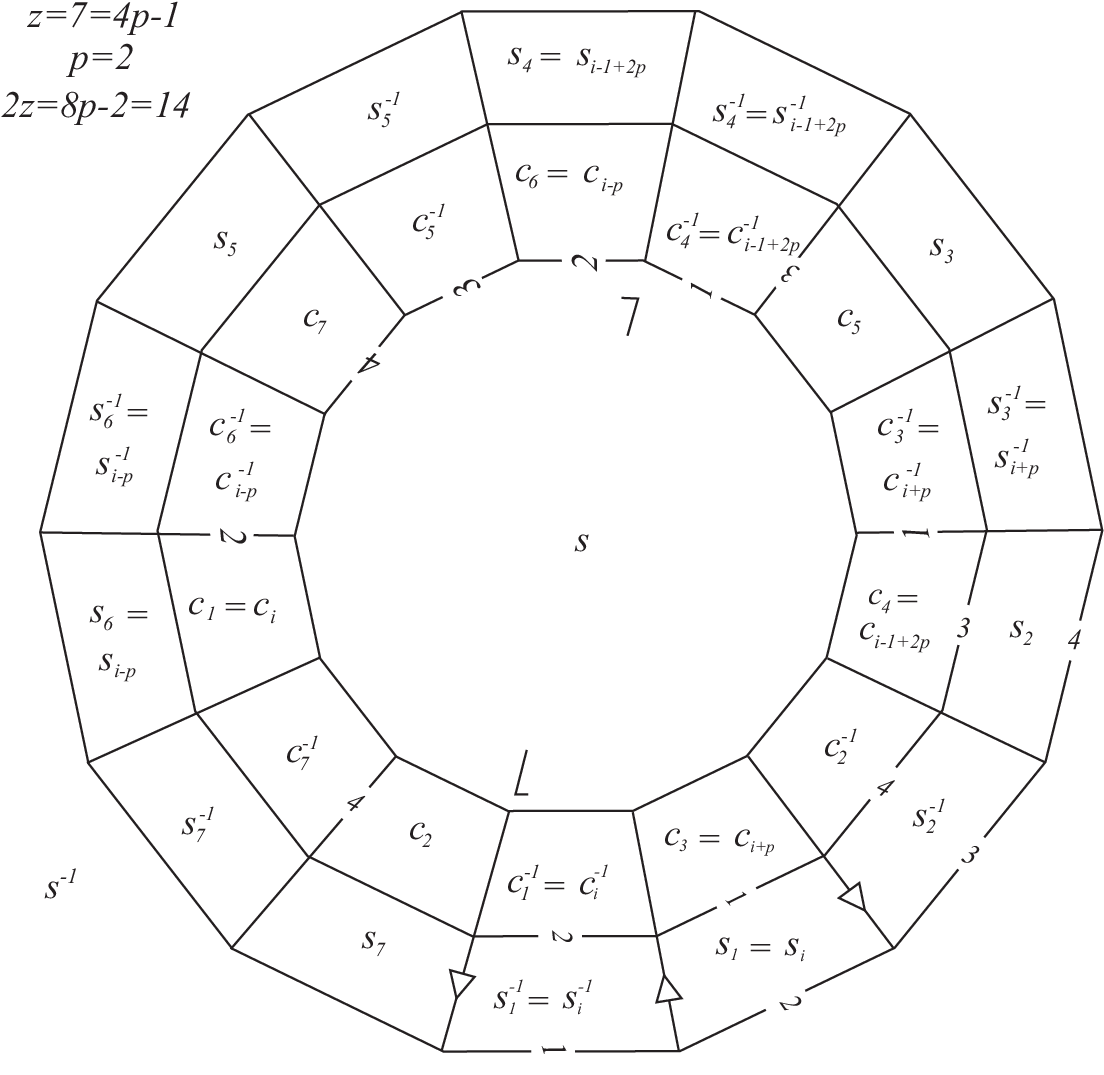}
\end{center}
\caption{The general new $Cw(2z)$ manifold scheme for $z = 4p-1$, illustrated by $p=2,z=7$.
See also formula (2.1) for algorithmic presentation}
\end{figure}

For instance the half screw
$\Bs_1: s_1^{-1} \rightarrow s_1$  for $\BCw(2z = 8p-2)$
with a fixed fundamental domain $W$ to faces $s_1^{-1}$ and $c_1^{-1}$
(see Fig.~4 and 6 together in ``mirror eyes") will be
$$\Bs_1 = \Ba_3(\Bh\Bb^1\Bh)\Bb^1\Bb^0$$
(see this also by the image of first edge $1$) on $s_1^{-1}$.
The screw motion $\Bc_1: c_1^{-1} \rightarrow c_1$ will
be
$$\Bc_1 = \Bh(\Bb^0 \Bb^1)^{1+p},$$
as you see by the edge $2$ on $s_1^{-1}$.
The tube screw
$\Bs: s^{-1} \rightarrow s$ will be
$$\Bs = \Ba_3(\Bh\Bb^0\Bh)(\Bb^1 \Bb^0)^{z-1}$$
by the image of edge $1$ on $s_1^{-1}$ again.
Our contribution (seems to be essential) to hyperbolic manifold theory is
\begin{figure}[]
\begin{center}
\includegraphics[width=13cm]{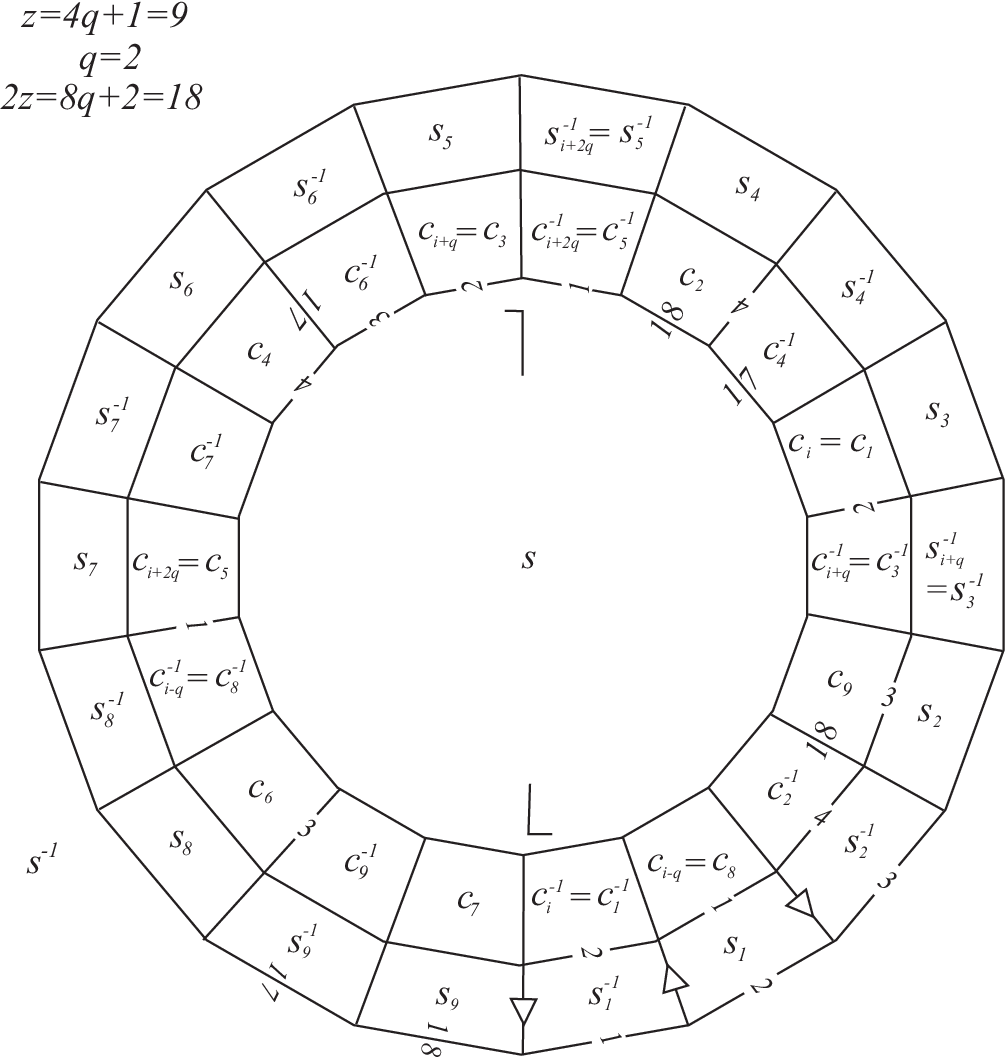}
\end{center}
\caption{General new manifold scheme for $z = 4q+1$, illustrated by $q=2$, i.e. $z=9$.
See also formula (2.2) for algorithmic presentation
}
\end{figure}
\begin{theorem}
The above cobweb (tube) manifold $Cw(2z)$
are minimal, i.e. none of them covers regularly another (smaller) manifold.
\end{theorem}
{\bf{Proof}} is a consequence of the Mostow rigidity theorem.
Namely, The fundamental group of a compact hyperbolic manifold
(as in our case as well) can be realized by a fixed-point-free isometry group of $\HYP$
with compact fundamental domain. The fundamental group of a regular covering manifold is an invariant subgroup of the fundamental group of the covered manifold. The factor group is called {\it covering group}.  Any symmetry map of $Cw(2z)$ is conjugated by an isometry that carries any orbit of $\BCw(2z)$, onto an orbit of it. So this holds for the generators of our half trunc-orthoscheme $W(2z)$, we can imagine it that tiles the fundamental domain $Cw(2z)$, as our manifold. That means, we would imagine two tilings of $Cw(2z)$: first, that with the fundamental domains of the covered manifolds under fixed-point-free isometries; second, with the half trunc-orthoscheme under plane reflections and halfturns. But this is not possible. The reflection domains cannot be divided into smaller parts, as well known for our orthoscheme. Other case cannot occur by the fixed point free actions of the covered manifold.
Contradiction!~ ~$\square$

\subsection{The strategy of our new dense ball packing construction}
Let $w=u$ and later $u=v=2z$, $JQEE_{13}F_{12}E_{02}F_{03} A_2$ be the vertices of the
half truncated orthoscheme $W(u;v;w=u)$ in Fig.~4.

All the essential typical point of $W$ can be expressed by its vector.
The types are the above vertices of $W$, a representing interior point of any edges,
included that of axis $h$, a representing point of any reflection face, a
representing point of the interior of $W$ can also be characterized by vectors, as we indicated in the Introduction (see also \cite{MSz17}, \cite{MSz18}). We do not repeat these here, since the computations will be implemented to computer.

E.g. point $A_2$ will be by the coordinate simplex vector $\BA_2$, and point $Q$
\begin{equation}
\begin{gathered}
Q(\BQ)=a_3 \cap A_3A_2; ~ \BQ=\BA_2-\frac{A_{23}}{A_{33}} \BA_3, \ \mathrm{with} \\
\langle \BQ,\BQ \rangle =\frac{(A_{22}A_{33}-A_{23}^2)}{A_{33}}=\langle \BQ,\BA_2 \rangle
=\frac{\sin^2\frac{\pi}{u}}{\sin^2\frac{\pi}{u}-\cos^2\frac{\pi}{v}}=\frac{A_{22}}{A_{33}}
\end{gathered} \notag
\end{equation}
by matrix (1.3).

Point $E(\BE)=a_3 \cap A_1A_3$ will be
\begin{equation}
\begin{gathered}
\BE=\BA_1-\frac{A_{13}}{A_{33}} \BA_3, \ \mathrm{with} \\
\langle \BE,\BE \rangle =\frac{(A_{11}A_{33}-A_{13}^2)}{A_{33}}=\langle \BE,\BA_1 \rangle
=\frac{1}{\sin^2\frac{\pi}{u}-\cos^2\frac{\pi}{v}}=\frac{1}{BA_{33}}.
\end{gathered} \notag
\end{equation}
In the considered cases $u=w$,
therefore, the midpoints $F_{03}$ of $A_0A_3$ and $F_{12}$ of $A_1A_2$,
respectively, can play important roles (?), since
$F_{03}F_{12}=h$ will be the axis of halfturn $\Bh$,
\begin{equation}
\begin{gathered}
\BF_{03}=\BA_0+\BA_3, ~ \langle \BF_{03},\BF_{03} \rangle= 2(A_{00}+A_{03})<0, \\
\BF_{12}=\BA_1+\BA_2, ~ \langle \BF_{12},\BF_{12} \rangle= 2(A_{11}+A_{12})<0.
\end{gathered} \notag
\end{equation}
So an interior point $H$ of axis $h$ can be $\BH = \BF_{03} + h \BF_{12}$
(with $0 < h$ real number). Also the other vertices, edge points, interior
points of reflection faces or interior points of the body $W$ can be described
by $0, 1, 2, 3$ real parameters, respectively. For their known distance formulas,
e.g. (1.4), (1.8) can be used for radii of packing balls, etc. For the contribution
$\delta(B_O)$ of an $O$-centred ball $B_O$ to packing density we shall use only
$W(u;v;w=u)$ by formula
\begin{equation}
\delta(B_O)=\frac{\mathrm{Vol}(B_O(r))}{|\BW_{O}|\mathrm{Vol}(W))}, \tag{2.4}
\end{equation}
where the maximal unified radius $r$ is restricted by the other balls,
$h$ and the walls of half trunc-orthoscheme $W$. The sum of the ball
contributions yield {\it the packing density in $W$ as our (new!?) definition
for packing with equal balls by different group orbits.}
For one ball orbit this agree with the general density definition of packing
not assumed to be regular but with equal balls (by D-V-cells, as in the introduction
of our paper \cite{MSz17}. But these are not valid for packings with equal balls
of more orbits (when the computation of the general density seems to be almost hopeless!?).
We recall some typical stabilizer orders (Fig.~1,4):
\begin{enumerate}
\item $A_2$ is the ball centre, $|\BW_{A_2}|=4u=|\BW_{Q}| ~ \text{for}~ O=Q$;
\item $F_{03}$ is the ball centre, $|\BW_{F_{03}}|=4v=|\BW_{J}| ~ \text{for}~ O=J$;
\item $F_{12}$ is the ball centre, $|\BW_{F_{12}}|=8=|\BW_{E}| ~ \text{for}~ O=E$.
\end{enumerate}
Of course, centre on $h$ or on reflection faces involves half ball,
centre in the interior of $W$ involves a full ball. It is easy to see,
that our manifold $Cw(2z)$, with $u= v = w = 2z$, $3 \le z$, by
its ``$2z$-gonal fundamental polyhedron" can be packed with full
equal balls, whose number of centre orbits is in accordance with
that in $W(2z)$, so that the above density in $W(2z)$ will be equal
to the sum of ball volumes divided ba the volume of $Cw(2z)$, just by symmetry argument.

{\bf{Strategies 2.2}}.
{\it Then the ball packing construction into a general but fixed $W(u; v;$ $w = u)$
for large density is a typical experimental interactive computer problem, e.g. by
two possible strategies:}
\begin{enumerate}
\item Say first one, as we follow: We compute the distances of vertices, edges,
faces from each other for overview. Then we place a centre into a vertex
(or edge point, face point) with a corresponding ball part of big as possible
radius, then try a full ball into the interior of big as possible radius.
Then equalize the two radii, and compute the packing density (by 2.4) with
two summands. Then we take two ball parts of big possible equal radii,
then try a full ball into the interior, etc. $\dots$, this is a finite
procedure. $W(2z)$ is extremely interesting for small $z$, say $3,5,7$.
\item Say second one, after distance overview: We place the first ball into
the interior of $W(u; v; w=u)$, then we take the ball parts,
radius equalization, densities, etc. ~ ~ $\square$
\end{enumerate}
This is a tedious work. A good program would be actual, and it seems to be important.
\section{Preliminary computation results}
We present here only our seemingly relevant results for later computations
by Strategy 1. Boldface
printed $\delta^{opt}$ density data and some less ones deserve further  approximations.!?
Thus, one ball part orbits seem to be not relevant,
compared with our former paper \cite{MSz17} adapted here.
\subsection{Two ball part orbits}
\medbreak
{\normalsize{
\centerline{\vbox{
\halign{\strut\vrule~\hfil $#$ \hfil~\vrule
&\quad \hfil $#$ \hfil~\vrule
&\quad \hfil $#$ \hfil\quad\vrule
&\quad \hfil $#$ \hfil\quad\vrule
&\quad \hfil $#$ \hfil\quad\vrule
\cr
\noalign{\hrule}
\noalign{\vskip2pt}
\multispan5{\strut\vrule\hfill\bf Table 1, Packings with 2 balls, $(u;v;w)=(3;7;3)$  \hfill\vrule}%
\cr
\noalign{\vskip2pt}
\noalign{\hrule}
\noalign{\vskip2pt}
\noalign{\hrule}
\text{Ball centres} & r^{opt} & 2 \cdot Vol(W_{uvw}) & Vol(B(r^{opt})) & \delta^{opt} \cr
\noalign{\hrule}
F_{12},E & 0.27140 & 0.27899 & 0.08498 & 0.15230 \cr
\noalign{\hrule}
F_{12},Q  & 0.27140 & 0.27899 & 0.08498 & 0.12692 \cr
\noalign{\hrule}
A_2,E  & 0.31215 & 0.27899 & 0.12991 & {\bf 0.19402} \cr
\noalign{\hrule}}}}
\smallbreak}}
\medbreak
{\normalsize{
\centerline{\vbox{
\halign{\strut\vrule~\hfil $#$ \hfil~\vrule
&\quad \hfil $#$ \hfil~\vrule
&\quad \hfil $#$ \hfil\quad\vrule
&\quad \hfil $#$ \hfil\quad\vrule
&\quad \hfil $#$ \hfil\quad\vrule
\cr
\noalign{\hrule}
\noalign{\vskip2pt}
\multispan5{\strut\vrule\hfill\bf Table 2, Packings with 2 balls, $(u;v;w)=(4;5;4)$  \hfill\vrule}%
\cr
\noalign{\vskip2pt}
\noalign{\hrule}
\noalign{\vskip2pt}
\noalign{\hrule}
Q,F_{03} & 0.53064 & 0.43062 & 0.66207 & 0.34594 \cr
\noalign{\hrule}
A_2,F_{03} & 0.53064 & 0.43062 & 0.66207 & 0.34594 \cr
\noalign{\hrule}
A_2,E  & 0.51921 & 0.43062 & 0.61872 & {\bf 0.53880} \cr
\noalign{\hrule}
A_2,J  & 0.53064 & 0.43062 & 0.66207 & 0.34594 \cr
\noalign{\hrule}}}}
\smallbreak}}
\medbreak
{\normalsize{
\centerline{\vbox{
\halign{\strut\vrule~\hfil $#$ \hfil~\vrule
&\quad \hfil $#$ \hfil~\vrule
&\quad \hfil $#$ \hfil\quad\vrule
&\quad \hfil $#$ \hfil\quad\vrule
&\quad \hfil $#$ \hfil\quad\vrule
\cr
\noalign{\hrule}
\noalign{\vskip2pt}
\multispan5{\strut\vrule\hfill\bf Table 3, Packings with 2 balls, $(u;v;w)=(5;4;5)$  \hfill\vrule}%
\cr
\noalign{\vskip2pt}
\noalign{\hrule}
\noalign{\vskip2pt}
\noalign{\hrule}
Q,F_{03} & 0.62687 & 0.46190 & 1.11606 & {\bf 0.54365} \cr
\noalign{\hrule}
A_2,F_{03} & 0.61123 & 0.46190 & 1.03061 & 0.50203 \cr
\noalign{\hrule}
F_{12},E & 0.41334 & 0.46190 & 0.30609 & 0.33134 \cr
\noalign{\hrule}
A_2,E  & 0.52717 & 0.46190 & 0.64869 & 0.49154 \cr
\noalign{\hrule}}}}
\smallbreak}}
\medbreak
{\normalsize{
\centerline{\vbox{
\halign{\strut\vrule~\hfil $#$ \hfil~\vrule
&\quad \hfil $#$ \hfil~\vrule
&\quad \hfil $#$ \hfil\quad\vrule
&\quad \hfil $#$ \hfil\quad\vrule
&\quad \hfil $#$ \hfil\quad\vrule
\cr
\noalign{\hrule}
\noalign{\vskip2pt}
\multispan5{\strut\vrule\hfill\bf Table 4, Packings with 2 balls, $(u;v;w)=(5;5;5)$  \hfill\vrule}%
\cr
\noalign{\vskip2pt}
\noalign{\hrule}
\noalign{\vskip2pt}
\noalign{\hrule}
Q,F_{03} & 0.56651 & 0.57271 & 0.81195 & 0.283558 \cr
\noalign{\hrule}
F_{12},E & 0.46453 & 0.57271 & 0.43838 & 0.38273 \cr
\noalign{\hrule}
A_2,E  & 0.60020 & 0.57271 & 0.97324 & {\bf 0.59477} \cr
\noalign{\hrule}}}}
\smallbreak}}
\medbreak
{\normalsize{
\centerline{\vbox{
\halign{\strut\vrule~\hfil $#$ \hfil~\vrule
&\quad \hfil $#$ \hfil~\vrule
&\quad \hfil $#$ \hfil\quad\vrule
&\quad \hfil $#$ \hfil\quad\vrule
&\quad \hfil $#$ \hfil\quad\vrule
\cr
\noalign{\hrule}
\noalign{\vskip2pt}
\multispan5{\strut\vrule\hfill\bf Table 5, Packings with 2 balls, $(u;v;w)=(6;4;6)$  \hfill\vrule}%
\cr
\noalign{\vskip2pt}
\noalign{\hrule}
\noalign{\vskip2pt}
\noalign{\hrule}
Q,F_{03} & 0.70337 & 0.55557& 1.60883 & 0.60329\cr
\noalign{\hrule}
A_2,F_{03} & 0.69217 & 0.55557& 1.52838 & 0.57313 \cr
\noalign{\hrule}
A_2,E  & 0.61947 & 0.55557 & 1.07504 & {\bf 0.64500} \cr
\noalign{\hrule}
A_2,J  & 0.65848 & 0.55557 & 1.30405 & 0.48901 \cr
\noalign{\hrule}}}}
\smallbreak}}
\medbreak
\medbreak
{\normalsize{
\centerline{\vbox{
\halign{\strut\vrule~\hfil $#$ \hfil~\vrule
&\quad \hfil $#$ \hfil~\vrule
&\quad \hfil $#$ \hfil\quad\vrule
&\quad \hfil $#$ \hfil\quad\vrule
&\quad \hfil $#$ \hfil\quad\vrule
\cr
\noalign{\hrule}
\noalign{\vskip2pt}
\multispan5{\strut\vrule\hfill\bf Table 6, Packings with 2 balls, $(u,v,w)=(7,4,7)$  \hfill\vrule}%
\cr
\noalign{\vskip2pt}
\noalign{\hrule}
\noalign{\vskip2pt}
\noalign{\hrule}
A_2,E  & 0.65278& 0.60917 & 1.26859 & {\bf 0.66938} \cr
\noalign{\hrule}}}}
\smallbreak}}
\medbreak
{\normalsize{
\centerline{\vbox{
\halign{\strut\vrule~\hfil $#$ \hfil~\vrule
&\quad \hfil $#$ \hfil~\vrule
&\quad \hfil $#$ \hfil\quad\vrule
&\quad \hfil $#$ \hfil\quad\vrule
&\quad \hfil $#$ \hfil\quad\vrule
\cr
\noalign{\hrule}
\noalign{\vskip2pt}
\multispan5{\strut\vrule\hfill\bf Table 7, Packings with 2 balls, $(u,v,w)=(7,3,7)$  \hfill\vrule}%
\cr
\noalign{\vskip2pt}
\noalign{\hrule}
\noalign{\vskip2pt}
\noalign{\hrule}
Q,F_{03} & 0.53278 & 0.38325 & 0.67042 & {\bf 0.41650} \cr
\noalign{\hrule}
A_2,F_{03} & 0.53278 & 0.38325 & 0.67042 & {\bf 0.41650} \cr
\noalign{\hrule}
F_{12},F_{03} & 0.36227 & 0.38325 & 0.20444 & 0.22227 \cr
\noalign{\hrule}
F_{12},E & 0.34463 & 0.38325 & 0.17557 & 0.22905 \cr
\noalign{\hrule}}}}
\smallbreak}}
\medbreak
{\normalsize{
\centerline{\vbox{
\halign{\strut\vrule~\hfil $#$ \hfil~\vrule
&\quad \hfil $#$ \hfil~\vrule
&\quad \hfil $#$ \hfil\quad\vrule
&\quad \hfil $#$ \hfil\quad\vrule
&\quad \hfil $#$ \hfil\quad\vrule
\cr
\noalign{\hrule}
\noalign{\vskip2pt}
\multispan5{\strut\vrule\hfill\bf Table 8, Packings with 2 balls, $(u,v,w)=(6,6,6)$  \hfill\vrule}%
\cr
\noalign{\vskip2pt}
\noalign{\hrule}
\noalign{\vskip2pt}
\noalign{\hrule}
F_{12},E & 0.46503 & 0.69130 & 0.43984 & 0.31812 \cr
\noalign{\hrule}
F_{12},Q  & 0.52985 & 0.69130 & 0.65901 & 0.31776 \cr
\noalign{\hrule}
A_2,E  & 0.57941 & 0.69130 & 0.87126 & {\bf 0.42010} \cr
\noalign{\hrule}}}}
\smallbreak}}
\medbreak
{\normalsize{
\centerline{\vbox{
\halign{\strut\vrule~\hfil $#$ \hfil~\vrule
&\quad \hfil $#$ \hfil~\vrule
&\quad \hfil $#$ \hfil\quad\vrule
&\quad \hfil $#$ \hfil\quad\vrule
&\quad \hfil $#$ \hfil\quad\vrule
\cr
\noalign{\hrule}
\noalign{\vskip2pt}
\multispan5{\strut\vrule\hfill\bf Table 9, Packings with 2 balls, $(u;v;w)=(10;10;10)$  \hfill\vrule}%
\cr
\noalign{\vskip2pt}
\noalign{\hrule}
\noalign{\vskip2pt}
\noalign{\hrule}
F_{12},E & 0.44294 & 0.83993 & 0.37858 & {\bf 0.22536} \cr
\noalign{\hrule}
F_{12},Q  & 0.32119 & 0.83993 & 0.14169 & 0.05061 \cr
\noalign{\hrule}
A_2,E  & 0.32119 & 0.83993 & 0.14169 & 0.05061 \cr
\noalign{\hrule}}}}
\smallbreak}}
\centerline{\vbox{
\halign{\strut\vrule~\hfil $#$ \hfil~\vrule
&\quad \hfil $#$ \hfil~\vrule
&\quad \hfil $#$ \hfil\quad\vrule
&\quad \hfil $#$ \hfil\quad\vrule
&\quad \hfil $#$ \hfil\quad\vrule
\cr
\noalign{\hrule}
\noalign{\vskip2pt}
\multispan5{\strut\vrule\hfill\bf Table 10, Packings with 2 balls, $(u;v;w)=(14;14;14)$  \hfill\vrule}%
\cr
\noalign{\vskip2pt}
\noalign{\hrule}
\noalign{\vskip2pt}
\noalign{\hrule}
F_{12},E & 0.44123 & 0.87770 & 0.37409 & {\bf 0.21311} \cr
\noalign{\hrule}
F_{12},Q  & 0.22661 & 0.87770 & 0.04925 & 0.01603 \cr
\noalign{\hrule}
A_2,E  & 0.22661 & 0.87770 & 0.04925 & 0.01603 \cr
\noalign{\hrule}}}}
\subsection{Three ball part orbits}
\medbreak
{\normalsize{
\centerline{\vbox{
\halign{\strut\vrule~\hfil $#$ \hfil~\vrule
&\quad \hfil $#$ \hfil~\vrule
&\quad \hfil $#$ \hfil\quad\vrule
&\quad \hfil $#$ \hfil\quad\vrule
&\quad \hfil $#$ \hfil\quad\vrule
\cr
\noalign{\hrule}
\noalign{\vskip2pt}
\multispan5{\strut\vrule\hfill\bf Table 11, Packings with 3 balls, $(u;v;w)=(3;7;3)$  \hfill\vrule}%
\cr
\noalign{\vskip2pt}
\noalign{\hrule}
\noalign{\vskip2pt}
\noalign{\hrule}
A_2,F_{03},E & 0.31215 & 0.27899 & 0.12991 & {\bf 0.22729} \cr
\noalign{\hrule}
Q,F_{12},F_{03} & 0.27140 & 0.27899 & 0.08498 & {0.14867} \cr
\noalign{\hrule}}}}
\smallbreak}}
\medbreak
{\normalsize{
\centerline{\vbox{
\halign{\strut\vrule~\hfil $#$ \hfil~\vrule
&\quad \hfil $#$ \hfil~\vrule
&\quad \hfil $#$ \hfil\quad\vrule
&\quad \hfil $#$ \hfil\quad\vrule
&\quad \hfil $#$ \hfil\quad\vrule
\cr
\noalign{\hrule}
\noalign{\vskip2pt}
\multispan5{\strut\vrule\hfill\bf Table 12, Packings with 3 balls, $(u;v;w)=(4;5;4)$  \hfill\vrule}%
\cr
\noalign{\vskip2pt}
\noalign{\hrule}
\noalign{\vskip2pt}
\noalign{\hrule}
A_2,F_{03},E & 0.51921 & 0.43062 & 0.61872 & {\bf 0.68248} \cr
\noalign{\hrule}
Q,F_{12},E & 0.38360 & 0.43062 & 0.24350 & 0.35341 \cr
\noalign{\hrule}
Q,F_{12},F_{03} & 0.38360 & 0.43062 & 0.24350 & 0.26859 \cr
\noalign{\hrule}}}}
\smallbreak}}
\medbreak
{\normalsize{
\centerline{\vbox{
\halign{\strut\vrule~\hfil $#$ \hfil~\vrule
&\quad \hfil $#$ \hfil~\vrule
&\quad \hfil $#$ \hfil\quad\vrule
&\quad \hfil $#$ \hfil\quad\vrule
&\quad \hfil $#$ \hfil\quad\vrule
\cr
\noalign{\hrule}
\noalign{\vskip2pt}
\multispan5{\strut\vrule\hfill\bf Table 13, Packings with 3 balls, $(u;v;w)=(5;4;5)$  \hfill\vrule}%
\cr
\noalign{\vskip2pt}
\noalign{\hrule}
\noalign{\vskip2pt}
\noalign{\hrule}
A_2,F_{03},E & 0.52717 & 0.46190 & 0.64869 & {\bf 0.66709} \cr
\noalign{\hrule}
A_2,F_{12},E  & 0.33587 & 0.46190 & 0.16233 &0.21087 \cr
\noalign{\hrule}
Q,F_{12},F_{03} & 0.41334 & 0.46190 & 0.30609 & 0.31477 \cr
\noalign{\hrule}}}}
\smallbreak}}
\medbreak
{\normalsize{
\centerline{\vbox{
\halign{\strut\vrule~\hfil $#$ \hfil~\vrule
&\quad \hfil $#$ \hfil~\vrule
&\quad \hfil $#$ \hfil\quad\vrule
&\quad \hfil $#$ \hfil\quad\vrule
&\quad \hfil $#$ \hfil\quad\vrule
\cr
\noalign{\hrule}
\noalign{\vskip2pt}
\multispan5{\strut\vrule\hfill\bf Table 14, Packings with 3 balls, $(u;v;w)=(5;5;5)$  \hfill\vrule}%
\cr
\noalign{\vskip2pt}
\noalign{\hrule}
\noalign{\vskip2pt}
\noalign{\hrule}
A_2,F_{03},E & 0.52477 & 0.57271 & 0.63955 & {\bf 0.50252} \cr
\noalign{\hrule}
A_2,F_{12},E  & 0.46453 & 0.57271 & 0.43838 & 0.45927 \cr
\noalign{\hrule}
Q,F_{12},F_{03} & 0.46453 & 0.57271 & 0.43838 & 0.34445 \cr
\noalign{\hrule}}}}
\smallbreak}}
\medbreak
{\normalsize{
\centerline{\vbox{
\halign{\strut\vrule~\hfil $#$ \hfil~\vrule
&\quad \hfil $#$ \hfil~\vrule
&\quad \hfil $#$ \hfil\quad\vrule
&\quad \hfil $#$ \hfil\quad\vrule
&\quad \hfil $#$ \hfil\quad\vrule
\cr
\noalign{\hrule}
\noalign{\vskip2pt}
\multispan5{\strut\vrule\hfill\bf Table 15, Packings with 3 balls, $(u;v;w)=(6;4;6)$  \hfill\vrule}%
\cr
\noalign{\vskip2pt}
\noalign{\hrule}
\noalign{\vskip2pt}
\noalign{\hrule}
A_2,F_{03},E & 0.49441 & 0.55557& 0.53155 & {\bf 0.43852} \cr
\noalign{\hrule}
Q,F_{12},E & 0.44069 & 0.55557& 0.37268 & 0.39130 \cr
\noalign{\hrule}
A_2,F_{12},E  & 0.42501 & 0.55557 & 0.33340 & 0.35006 \cr
\noalign{\hrule}
Q,F_{12},F_{03} & 0.46151 & 0.55557 & 0.42966 & 0.35446\cr
\noalign{\hrule}}}}
\smallbreak}}
\medbreak
\medbreak
{\normalsize{
\centerline{\vbox{
\halign{\strut\vrule~\hfil $#$ \hfil~\vrule
&\quad \hfil $#$ \hfil~\vrule
&\quad \hfil $#$ \hfil\quad\vrule
&\quad \hfil $#$ \hfil\quad\vrule
&\quad \hfil $#$ \hfil\quad\vrule
\cr
\noalign{\hrule}
\noalign{\vskip2pt}
\multispan5{\strut\vrule\hfill\bf Table 16, Packings with 3 balls, $(u;v;w)=(7;4;7)$  \hfill\vrule}%
\cr
\noalign{\vskip2pt}
\noalign{\hrule}
\noalign{\vskip2pt}
\noalign{\hrule}
Q,F_{12},E  & 0.48955 & 0.60917 & 0.51555 & {\bf 0.48361} \cr
\noalign{\hrule}
A_2,F_{12},E  & 0.48955 & 0.60917 & 0.51555 & {\bf 0.48361} \cr
\noalign{\hrule}
Q,F_{12},F_{03} & 0.48955 & 0.60917 & 0.51555 &  0.37782 \cr
\noalign{\hrule}}}}
\smallbreak}}
\medbreak
{\normalsize{
\centerline{\vbox{
\halign{\strut\vrule~\hfil $#$ \hfil~\vrule
&\quad \hfil $#$ \hfil~\vrule
&\quad \hfil $#$ \hfil\quad\vrule
&\quad \hfil $#$ \hfil\quad\vrule
&\quad \hfil $#$ \hfil\quad\vrule
\cr
\noalign{\hrule}
\noalign{\vskip2pt}
\multispan5{\strut\vrule\hfill\bf Table 17, Packings with 3 balls, $(u;v;w)=(7;3;7)$  \hfill\vrule}%
\cr
\noalign{\vskip2pt}
\noalign{\hrule}
\noalign{\vskip2pt}
\noalign{\hrule}
A_2,F_{12},E & 0.38521 & 0.38325 & 0.24664 & {\bf 0.36774} \cr
\noalign{\hrule}
Q,F_{12},F_{03} & 0.36227 & 0.38325 & 0.20444 & 0.26038 \cr
\noalign{\hrule}}}}
\smallbreak}}
\medbreak
{\normalsize{
\centerline{\vbox{
\halign{\strut\vrule~\hfil $#$ \hfil~\vrule
&\quad \hfil $#$ \hfil~\vrule
&\quad \hfil $#$ \hfil\quad\vrule
&\quad \hfil $#$ \hfil\quad\vrule
&\quad \hfil $#$ \hfil\quad\vrule
\cr
\noalign{\hrule}
\noalign{\vskip2pt}
\multispan5{\strut\vrule\hfill\bf Table 18, Packings with 3 balls, $(u;v;w)=(6;6;6)$  \hfill\vrule}%
\cr
\noalign{\vskip2pt}
\noalign{\hrule}
\noalign{\vskip2pt}
\noalign{\hrule}
Q,F_{12},E & 0.46503 & 0.69130 & 0.43984 & {\bf 0.37114} \cr
\noalign{\hrule}
A_2,F_{12},E & 0.46503 & 0.69130 & 0.43984 & {\bf 0.37114} \cr
\noalign{\hrule}}}}
\smallbreak}}
\medbreak
{\normalsize{
\centerline{\vbox{
\halign{\strut\vrule~\hfil $#$ \hfil~\vrule
&\quad \hfil $#$ \hfil~\vrule
&\quad \hfil $#$ \hfil\quad\vrule
&\quad \hfil $#$ \hfil\quad\vrule
&\quad \hfil $#$ \hfil\quad\vrule
\cr
\noalign{\hrule}
\noalign{\vskip2pt}
\multispan5{\strut\vrule\hfill\bf Table 19, Packings with 3 balls, $(u;v;w)=(10;10;10)$  \hfill\vrule}%
\cr
\noalign{\vskip2pt}
\noalign{\hrule}
\noalign{\vskip2pt}
\noalign{\hrule}
Q,F_{12},E & 0.32119 & 0.83993 & 0.14169 & {\bf 0.09278} \cr
\noalign{\hrule}
A_2,F_{12},E & 0.32119 & 0.83993 & 0.14169 & {\bf 0.09278} \cr
\noalign{\hrule}}}}
\smallbreak}}
\medbreak
{\normalsize{
\centerline{\vbox{
\halign{\strut\vrule~\hfil $#$ \hfil~\vrule
&\quad \hfil $#$ \hfil~\vrule
&\quad \hfil $#$ \hfil\quad\vrule
&\quad \hfil $#$ \hfil\quad\vrule
&\quad \hfil $#$ \hfil\quad\vrule
\cr
\noalign{\hrule}
\noalign{\vskip2pt}
\multispan5{\strut\vrule\hfill\bf Table 20, Packings with 3 balls, $(u;v;w)=(14;14;14)$  \hfill\vrule}%
\cr
\noalign{\vskip2pt}
\noalign{\hrule}
\noalign{\vskip2pt}
\noalign{\hrule}
Q,F_{12},E & 0.22661 & 0.87770 & 0.04925 & {\bf 0.03006} \cr
\noalign{\hrule}
A_2,F_{12},E & 0.22661 & 0.87770 & 0.04925 & {\bf 0.03006} \cr
\noalign{\hrule}}}}
\smallbreak}}
\subsection{Four ball part orbits}
\medbreak
{\normalsize{
\centerline{\vbox{
\halign{\strut\vrule~\hfil $#$ \hfil~\vrule
&\quad \hfil $#$ \hfil~\vrule
&\quad \hfil $#$ \hfil\quad\vrule
&\quad \hfil $#$ \hfil\quad\vrule
&\quad \hfil $#$ \hfil\quad\vrule
\cr
\noalign{\hrule}
\noalign{\vskip2pt}
\multispan5{\strut\vrule\hfill\bf Table 21, Packings with 4 balls, $(u;v;w)=(3;7;3)$  \hfill\vrule}%
\cr
\noalign{\vskip2pt}
\noalign{\hrule}
\noalign{\vskip2pt}
\noalign{\hrule}
A_2,F_{03},E,F_{12} & 0.15608 & 0.27899 & 0.01600 & {\bf 0.04234} \cr
\noalign{\hrule}}}}
\smallbreak}}
\medbreak
{\normalsize{
\centerline{\vbox{
\halign{\strut\vrule~\hfil $#$ \hfil~\vrule
&\quad \hfil $#$ \hfil~\vrule
&\quad \hfil $#$ \hfil\quad\vrule
&\quad \hfil $#$ \hfil\quad\vrule
&\quad \hfil $#$ \hfil\quad\vrule
\cr
\noalign{\hrule}
\noalign{\vskip2pt}
\multispan5{\strut\vrule\hfill\bf Table 22, Packings with 4 balls, $(u;v;w)=(4;5;4)$  \hfill\vrule}%
\cr
\noalign{\vskip2pt}
\noalign{\hrule}
\noalign{\vskip2pt}
\noalign{\hrule}
A_2,F_{03},E,F_{12} & 0.26532 & 0.43062 & 0.07934 & {\bf 0.13358} \cr
\noalign{\hrule}
Q,F_{12},E,F_{03} & 0.26532 & 0.43062 & 0.07934 & {\bf 0.13358} \cr
\noalign{\hrule}}}}
\smallbreak}}
\medbreak
{\normalsize{
\centerline{\vbox{
\halign{\strut\vrule~\hfil $#$ \hfil~\vrule
&\quad \hfil $#$ \hfil~\vrule
&\quad \hfil $#$ \hfil\quad\vrule
&\quad \hfil $#$ \hfil\quad\vrule
&\quad \hfil $#$ \hfil\quad\vrule
\cr
\noalign{\hrule}
\noalign{\vskip2pt}
\multispan5{\strut\vrule\hfill\bf Table 23, Packings with 4 balls, $(u;v;w)=(5;4;5)$  \hfill\vrule}%
\cr
\noalign{\vskip2pt}
\noalign{\hrule}
\noalign{\vskip2pt}
\noalign{\hrule}
A_2,F_{03},E,F_{12} & 0.33587 & 0.46190 & 0.16233 & {\bf 0.25480} \cr
\noalign{\hrule}
Q,F_{12},E,F_{03} & 0.31344 & 0.46190 & 0.13154 & 0.20647 \cr
\noalign{\hrule}}}}
\smallbreak}}

\medbreak
{\normalsize{
\centerline{\vbox{
\halign{\strut\vrule~\hfil $#$ \hfil~\vrule
&\quad \hfil $#$ \hfil~\vrule
&\quad \hfil $#$ \hfil\quad\vrule
&\quad \hfil $#$ \hfil\quad\vrule
&\quad \hfil $#$ \hfil\quad\vrule
\cr
\noalign{\hrule}
\noalign{\vskip2pt}
\multispan5{\strut\vrule\hfill\bf Table 24, Packings with 4 balls, $(u;v;w)=(5;5;5)$  \hfill\vrule}%
\cr
\noalign{\vskip2pt}
\noalign{\hrule}
\noalign{\vskip2pt}
\noalign{\hrule}
A_2,F_{03},E,F_{12} & 0.37373 & 0.57271 & 0.22485 & {0.27483} \cr
\noalign{\hrule}
Q,F_{12},E,F_{03} & 0.42124 & 0.57271 & 0.32440 & {\bf 0.39650} \cr
\noalign{\hrule}}}}
\smallbreak}}
\medbreak
{\normalsize{
\centerline{\vbox{
\halign{\strut\vrule~\hfil $#$ \hfil~\vrule
&\quad \hfil $#$ \hfil~\vrule
&\quad \hfil $#$ \hfil\quad\vrule
&\quad \hfil $#$ \hfil\quad\vrule
&\quad \hfil $#$ \hfil\quad\vrule
\cr
\noalign{\hrule}
\noalign{\vskip2pt}
\multispan5{\strut\vrule\hfill\bf Table 25, Packings with 4 balls, $(u;v;w)=(6;4;6)$  \hfill\vrule}%
\cr
\noalign{\vskip2pt}
\noalign{\hrule}
\noalign{\vskip2pt}
\noalign{\hrule}
A_2,F_{03},E,F_{12} & 0.42501 & 0.55557 & 0.33340 & 0.42507 \cr
\noalign{\hrule}
Q,F_{12},E,F_{03} & 0.44069 & 0.55557 & 0.37268 & {\bf 0.47515} \cr
\noalign{\hrule}}}}
\smallbreak}}
\medbreak
{\normalsize{
\centerline{\vbox{
\halign{\strut\vrule~\hfil $#$ \hfil~\vrule
&\quad \hfil $#$ \hfil~\vrule
&\quad \hfil $#$ \hfil\quad\vrule
&\quad \hfil $#$ \hfil\quad\vrule
&\quad \hfil $#$ \hfil\quad\vrule
\cr
\noalign{\hrule}
\noalign{\vskip2pt}
\multispan5{\strut\vrule\hfill\bf Table 26, Packings with 4 balls, $(u;v;w)=(7;4;7)$  \hfill\vrule}%
\cr
\noalign{\vskip2pt}
\noalign{\hrule}
\noalign{\vskip2pt}
\noalign{\hrule}
A_2,F_{03},E,F_{12} & 0.48947 & 0.60917 & 0.51528 & {\bf 0.58910} \cr
\noalign{\hrule}
Q,F_{12},E,F_{03} & 0.48947 & 0.60917 & 0.51528 & {\bf 0.58910} \cr
\noalign{\hrule}}}}
\smallbreak}}
\medbreak
{\normalsize{
\centerline{\vbox{
\halign{\strut\vrule~\hfil $#$ \hfil~\vrule
&\quad \hfil $#$ \hfil~\vrule
&\quad \hfil $#$ \hfil\quad\vrule
&\quad \hfil $#$ \hfil\quad\vrule
&\quad \hfil $#$ \hfil\quad\vrule
\cr
\noalign{\hrule}
\noalign{\vskip2pt}
\multispan5{\strut\vrule\hfill\bf Table 27, Packings with 4 balls, $(u;v;w)=(7;3;7)$  \hfill\vrule}%
\cr
\noalign{\vskip2pt}
\noalign{\hrule}
\noalign{\vskip2pt}
\noalign{\hrule}
A_2,F_{03},E,F_{12} & 0.34463 & 0.38325 & 0.17557 & {\bf 0.33812} \cr
\noalign{\hrule}
Q,F_{12},E,F_{03} & 0.34463 & 0.38325 &0.17557 & {\bf 0.33812} \cr
\noalign{\hrule}
J,Q,F_{12},A_2 & 0.28313 & 0.38325 & 0.09661 & {0.14104} \cr
\noalign{\hrule}}}}
\smallbreak}}
\medbreak
{\normalsize{
\centerline{\vbox{
\halign{\strut\vrule~\hfil $#$ \hfil~\vrule
&\quad \hfil $#$ \hfil~\vrule
&\quad \hfil $#$ \hfil\quad\vrule
&\quad \hfil $#$ \hfil\quad\vrule
&\quad \hfil $#$ \hfil\quad\vrule
\cr
\noalign{\hrule}
\noalign{\vskip2pt}
\multispan5{\strut\vrule\hfill\bf Table 28, Packings with 4 balls, $(u;v;w)=(6;6;6)$  \hfill\vrule}%
\cr
\noalign{\vskip2pt}
\noalign{\hrule}
\noalign{\vskip2pt}
\noalign{\hrule}
A_2,F_{03},E,F_{12} & 0.37764 & 0.691309 & 0.23211 & {\bf 0.22384} \cr
\noalign{\hrule}
Q,F_{12},E,F_{03} & 0.37764 & 0.691309 & 0.23211 & {\bf 0.22384} \cr
\noalign{\hrule}}}}
\smallbreak}}
\medbreak
{\normalsize{
\centerline{\vbox{
\halign{\strut\vrule~\hfil $#$ \hfil~\vrule
&\quad \hfil $#$ \hfil~\vrule
&\quad \hfil $#$ \hfil\quad\vrule
&\quad \hfil $#$ \hfil\quad\vrule
&\quad \hfil $#$ \hfil\quad\vrule
\cr
\noalign{\hrule}
\noalign{\vskip2pt}
\multispan5{\strut\vrule\hfill\bf Table 29, Packings with 4 balls, $(u;v;w)=(10;10;10)$  \hfill\vrule}%
\cr
\noalign{\vskip2pt}
\noalign{\hrule}
\noalign{\vskip2pt}
\noalign{\hrule}
A_2,F_{03},E,F_{12} & 0.17700 & 0.83993 & 0.023373 & {\bf 0.01670} \cr
\noalign{\hrule}
Q,F_{12},E,F_{03} & 0.17700 & 0.83993 & 0.023373 & {\bf 0.01670} \cr
\noalign{\hrule}}}}
\smallbreak}}
\medbreak
{\normalsize{
\centerline{\vbox{
\halign{\strut\vrule~\hfil $#$ \hfil~\vrule
&\quad \hfil $#$ \hfil~\vrule
&\quad \hfil $#$ \hfil\quad\vrule
&\quad \hfil $#$ \hfil\quad\vrule
&\quad \hfil $#$ \hfil\quad\vrule
\cr
\noalign{\hrule}
\noalign{\vskip2pt}
\multispan5{\strut\vrule\hfill\bf Table 30, Packings with 4 balls, $(u;v;w)=(14;14;14)$  \hfill\vrule}%
\cr
\noalign{\vskip2pt}
\noalign{\hrule}
\noalign{\vskip2pt}
\noalign{\hrule}
A_2,F_{03},E,F_{12} & 0.11911 & 0.87770 & 0.00710 & {\bf 0.00462} \cr
\noalign{\hrule}
Q,F_{12},E,F_{03} & 0.11911 & 0.87770 & 0.00710 & {\bf 0.00462}  \cr
\noalign{\hrule}}}}
\smallbreak}}
\section{Summary}
We briefly collect our new main results in two summarizing theorems

\begin{theorem}
The cobweb (tube) manifolds $Cw(2z)$ to fundamental cobweb (tube)
polyhedra as fundamental domains have been constructed by two series of simplified
face pairing identifications in Figures $4, 5, 6, 7,$ described above in Section 2.
Any fundamental group $\BCw(2z = 8p - 2)$ of first series and $\BCw(2z = 8q + 2)$
of second series can be given by the algorithmic presentation in formulas
(2.1) and (2.2), to Figures $6$ and $7$, respectively, on the base of their
symmetry group $\BW(2z)$ with half trunc-orthoscheme fundamental domain $W(2z)$ in
Fig.~4.

All necessary metric data of $Cw(2z)$ can be computed on the base of
this half trunc-orthoscheme and its projectiv-metric coordinate simplex
$b^0b^1b^2b^3 = A_0 A_1 A_2A_3$ by scalar products in (1.1), (1.3),
respectively, or their (polar) plane $\rightarrow$ (pole) point polarity in Section 1.
~ ~ ~ $\square$
\end{theorem}
\begin{theorem}
The above cobweb (tube) manifolds $Cw(2z)$ are minimal, i.e. none of them covers
regularly another (smaller) manifold. This assertion is based on the symmetry group
$\BW(2z)$ that is an extended complete Coxeter reflection group with the half
trunc-orthoscheme $W(2z = u = v = w)$ as fundamental domain in Figures 1, 4.
So we also recall a most economical presentation of the general extended complete
reflection group $\BW(u; v; w=u)$ in (2.3).

We have given strategies in Section 2.2 for dense ball packing constructions
with equal balls belonging to more ball centre orbits by $\BCw(2z)$, or in more
general by $\BW(u; v; w=u)$. As first results, we give dense enough constructions
by the first strategy in Section 3. ~ ~ ~ $\square$
\end{theorem}
For instance, in Table 12 at $(u; v; w =u) = (4; 5; 4)$ we look our top optimum $\delta^{opt} \approx 0.68248$
for the $3$ orbits of $A_2, F_{03}, E$ in the half trunc-orthoscheme $W$ with $1/16, 1/20, 1/8$
ball parts, respectively of optimal radius $\approx 0.51921$, etc, that yield this preliminary optimum.
This can be increased if an apropriate full ball can be placed into an interior pont $I$ of $W$.

At Table $8$ of $(6; 6; 6)$ we similarly look $\delta^{opt} \approx 0.42010$ for $2$ orbits of
$A_2, E$ in $W$
with $1/24, 1/8$ ball parts, respectively with optimal radius $\approx 0.57941$, etc.
Thus we get a preliminary optimum, that increases if we would place an appropriate full ball into $W$.
Then we can exactly tell, how many equal balls are placed into the above tube manifold $Cw(6)$.

We intend to continue these investigations, also for promising probable material applications.

{\bf We would like to Congratulate the Staff of Matemati\v cki Vesnik on the nice
75th Anniversary with our Sincere Best Wishes!}


\begin{thebibliography}{999999999}
%
\bibitem{CaTe09}
{Cavicchioli,~A.~--~ Telloni,~ A. I.,}
On football manifolds of E. Moln\'ar.
\textit{Acta Math.Hungar.,}
{\bf 124(4)} (2009), 321-332.
%
\bibitem{K89} Kellerhals,~R. {On the volume of hyperbolic polyhedra},
{\em Math. Ann.}, {\bf 245} (1989), \rm 541--569.
%
\bibitem{M89} Moln\'ar,~E., {Projective metrics and hyperbolic volume,}
\textit{Annales Univ. Sci. Budapest, Sect. Math.,},  {\bf{32}}, (1989), p. 127-157.
%
\bibitem{M84}
 {Moln\'ar,~E.,}
 Space forms and fundamental polyhedra.
 \textit{ Proceedings of the Conference on Differential
 Geometry and Its Applications, Nov\'e M\'esto na Morav\'e, Czechoslovakia
 1983. Part 1. Differential Geometry,}  (1984), 91-103.
%
\bibitem{M88}
{Moln\'ar,~E.,}
Two hyperbolic football manifolds.
\textit{ Proceedings of International Conference on Differential
Geometry and Its Applications, Dubrovnik Yugoslavia,} (1988), 217-241.
%
\bibitem{M05} Moln\'ar,~E.,{Combinatorial construction of tilings by barycentric simplex orbits
(D-symbols) and their realizations in Euclidean and other homogeneous spaces}.
\textit{Acta Cryst.}, (2005) \bf{A61}\rm, 541--552.
%
\bibitem{M97} Moln\'ar,~E. The projective interpretation of the eight 3-dimensional
homogeneous geometries,
{\em{Beitr. Alg. Geom., (Contr. Alg. Geom.)}}, {\bf 38/2} (1997), 261--288.
%
\bibitem{M92}
 {Moln\'ar,~E.,}
 Polyhedron complexes with simply transitive group actions and their realizations.
 \textit{Acta Math. Hung.,}  {\bf 59 (1-2)} (1992), 175-216.
%
\bibitem{M91}
{Moln\'ar,~E.,}
{On non-Euclidean crystallography, some football manifolds,}
\textit{Structural Chemistry,} {\bf{23/4}}, (2012), 1057-1069.
%
\bibitem{MSz}
{Moln{\'a}r,~E.~--~Szirmai,~J.,}
Symmetries in the 8 homogeneous 3-geometries.
\textit{Symmetry Cult. Sci.,}
{\bf 21/1-3} (2010), 87-117.
%
\bibitem{MSz17}
{Moln{\'a}r,~E.~--~Szirmai,~J.,}
Top dense hyperbolic ball packings and coverings for complete Coxeter
orthoscheme groups.
\textit{Publications de l'Institut Mathématique,}
Nouvelle série, tome {\bf 103 (117)} (2018), 129-146.
%
\bibitem{MSz16}
{Moln{\'a}r,~E.~--~Szirmai,~J.,}
On hyperbolic cobweb manifolds.
{\em {Stud. Univ. Zilina. Math .Ser.}}, {\bf 28}, (2016), 43--52.
%
\bibitem{MSz18}
{Moln{\'a}r,~E.~--~Szirmai,~J.,}
Infinite series of compact hyperbolic manifolds,
as possible crystal structures.
{\em {Matemati\v cki Vesnik}}, {bf 72, 3}, (2020), 257-272.
%
\bibitem{M-Sz}
{Moln{\'a}r,~E.~--~Szirmai,~J.,}
On homogeneous 3-geometries, balls and their optimal arrangements, especially
in $\NIL$ and $\SOL$ spaces.
\emph{G-Slovak Journal for Geometry and Graphics} {\bf 19}  (2022), 37, 5-32.
%
\bibitem{stachel}
Moln{\'a}r,~E.~--~Szirmai,~J.,
Packings with geodesic and translation balls and their visualizations in $\SLR$ space.
{\it Journal for Geometry and Graphics}  {\bf 26}(1) (2022), 51--64.
%
\bibitem{Boh}
Moln{\'a}r,~E.~--~Szirmai,~J.,
Non-Euclidean Crystal Geometry,
{\it South Bohemia Mathematical Letters} {\bf 30}, No.1.  (2022), 28--40.
%
\bibitem{Mong23}
Moln{\'a}r,~E.~---~Szirmai,~J.,
On $\SLR$ crystallography.
{\it Proceedings  of the 9th International Scientific Conference
on Geometry and Graphics MoNGeometrija 2023,
Geometry, graphics and design in the digital age, I. Bajsanski,
M. Jovanovic [ed],
Faculty of Technical Sciences, University of Novi Sad
Serbian Society for Geometry and Graphics SUGIG} (2023), 229--245.
%
\bibitem{MPSz06}
Moln{\'a}r,~E.~--~Prok,~I.~--~Szirmai,~J. Classification of
tile-transitive 3-simplex tilings and their realizations in
homogeneous spaces, {\em {A.~Pr\'ekopa and E.~Moln\'ar, (eds.).
Non-Euclidean Geometries, J\'anos Bolyai Memorial Volume,
Mathematics and Its Applications}}, Springer, (2006), Vol.~{\bf 581}, 321--363.
%
\bibitem{MStSz}
Moln{\'a}r,~E.~--~Stojanovi\'c,~M.~--~Szirmai,~J.
Non-fundamental trunc-simplex tilings and their optimal hyperball packings and coverings
in hyperbolic space \\ I. For families F1-F4,
\emph{Filomat}, {\bf37:5} (2023), 1409--1448.
%
\bibitem{P92} Prok,~I. Data structures and procedures for a polyhedron algorithm,
{\em {Periodica Polytechnica Ser. Mech. Eng.}}, {\bf 36/3-4}, (1992), 299--316.
%
\bibitem{P98} Prok,~I. Classification of dodecahedral space forms,
{\em {Beitr. Alg. Geom., (Contr. Alg. Geom.)}}, {\bf 38/2} (1998), 497--515.
%
\bibitem{S83}
{Scott,~P.}
The geometries of 3-manifolds.
{\em{Bull. London Math. Soc.}}, {\bf 15} (1983), 401--487.
%
\bibitem{Sz07-2} Szirmai,~J. The densest geodesic ball packing by a type of $\NIL$ lattices,
{\em{Beitr. Alg. Geom., (Contr. Alg. Geom.)}}, {\bf 48/2}, (2007), 383--397.
%
\bibitem{Sz13-2}
{Szirmai,~J.}
 A candidate to the densest packing with equal balls in the Thurston geometries.
{\em {Beitr. Algebra Geom.,}} {\bf 55/2}, (2014), 441- 452.
%
\bibitem{V93}
{Vinberg,~E.~B. (Ed.)}
\textit{Geometry II. Spaces of Constant Curvature}
Spriger Verlag Berlin-Heidelberg-New York-London-Paris-Tokyo-Hong Kong-Barcelona-Budapest, 1993.
%
\bibitem{W06}
Weeks,~J.~R.
{Real-time animation in hyperbolic, spherical, and product geometries.}
{\em {A.~Pr\'ekopa and E.~Moln\'ar, (eds.).
Non-Euclidean Geometries, J\'anos Bolyai Memorial Volume,
Mathematics and Its Applications}}, Springer (2006) Vol.~{\bf 581}, 287--305.
%
\bibitem{W}
 {Wolf,~J.~A.}
\textit{Spaces of Constant Curvature.}
 McGraw-Hill, New York, 1967, (Russian translation: Izd. "Nauka" Moscow, 1982).
%
\end{thebibliography}
\end{document}